\documentclass[a4paper]{amsart}
\usepackage{amsmath}
\def\h{\mathcal{H}}
\def\s3{\mathcal{S}_3}
\def\dimo{\textit{Proof}$\quad$}

\def\Torus{\mathbb{T}}

\def\contr{\neg}
\def\dime{3}

\def\contr{\rightharpoonup}

\def\R{\mathbb{R}}
\def\C{\mathbb{C}}
\def\Z{\mathbb{Z}}

\def\gg{\mathcal{L}_\C}
\def\ggr{\mathcal{L}}
\def\so2{\mathbf{so}(2,\R)}
\def\sl2c{\mathbf{sl}(2,\C)}
\def\slsc{\mathbf{sl}(6,\C)}
\def\ka{K\"{a}hler\;}

\newtheorem{teo}{Theorem}[section]

\newtheorem{cor}[teo]{Corollary}

\newtheorem{lem}[teo]{Lemma}
\newtheorem{pro}[teo]{Proposition}
\newtheorem{dfn}[teo]{Definition}

\newtheorem{rmk}[teo]{Remark}

\title{A geometric realization of $\slsc$}
\author{Giovanni Gaiffi, ~Michele Grassi}
\date{October 24, 2006}
\begin{document}
\begin{abstract}
Given an orientable weakly self-dual manifold $X$ of rank two, we
build a geometric realization of the Lie algebra $\slsc$ as a naturally defined
algebra $\gg$ of endomorphisms of the space of differential forms of
$X$. We provide an explicit description of  Serre generators in
terms of natural generators of $\gg$. This construction gives a
bundle on $X$ which is related to the search for a natural Gauge
theory on $X$. We consider this paper as a first step in the study of a rich and interesting algebraic
structure.
\end{abstract}
\maketitle
\section{Introduction}
This paper is a step in a broader program, which aims at finding a
geometric counterpart to the Mirror Symmetry phaenomenon, and
possibly a geometric language in which to formulate a physical
theory interpolating between different $\sigma$-models. While we
direct the reader to \cite{G2},\cite{G3} for more details, we list
here only some aspects of this theory to put the present work into
context.

In the Strominger-Yau-Zaslow approach to Mirror Symmetry you have
that two mirror dual Calabi-Yaus should posses (in some limiting
sense) semi-flat special lagrangian torus fibrations $f:M\to B$,
$\hat{f}:\hat{M}\to B$ which have as fibres flat tori which are dual
in the metric sense (see \cite{SYZ}, and \cite{G2} for the
terminology and the definitions). As it is widely known, the major
drawback of this approach is that it is very difficult to build
special lagrangian tori fibrations. Usually this construction can be
carried out only when the dual Calabi-Yau manifolds are actually
hyperkahler, and the special lagrangian tori can be viewed as
complex submanifolds (with respect to a rotated complex structure),
so that
the methods of complex algebraic geometry can be put to work.\\
When you do have the fibrations, then the idea is to construct the
mirror map as a sort of Fourier-Mukai transform (see for example
\cite{BMP}). This Fourier-Mukai transform is a correspondence
induced by pull-back and push forward from the space $X = M\times_B
\hat{M}$. In the hyperk\"ahler case this space is a complex
manifold, while in the general case (for example for Mirror Symmetry
for Calabi-Yau threefolds) it is just a real manifold of (real)
dimension $3\cdot dim_\C(M)$.\\

\textbf{Background.} The notion of (Weakly) self-dual manifold (cf. \cite{G2}) was
conceived in the first place to isolate the geometric aspects of the
$X$ above which are needed to obtain Mirror Symmetry between $M$ and
$\hat{M}$. We reproduce here the definition for the reader, while
referring to \cite{G2} and \cite{G3} for all the remarks, examples and
observations:
\begin{dfn}
A {\em weakly self-dual manifold} (WSD manifold for brevity) is
given by a smooth manifold $X$, together with two smooth $2$-forms
$\omega_1,\omega_2$ a Riemannian metric and a third smooth
$2$-form $\omega_D$ (the {\em dualizing} form) on it, which satisfy the following conditions:\\
1) $d\omega_1 = d\omega_2 = d\omega_D = 0$ and the distribution $\omega_1^0+\omega_2^0$ is integrable.\\
2) For all $p\in X$ there exist an orthogonal basis
$dx_1,..,dx_m,dy^1_1,...,dy^1_m,dy^2_1,...,dy^2_m$,
$dz_1,...,dz_c,dw_1,...,dw_c$ of  $T_p^*X$ such that the
$dx_1,..,dx_m,dy^1_1,...,dy^1_m,dy^2_1,...,dy^2_m$ are orthonormal
and
\[(\omega_1)_p = \sum_{i=1}^m dx_i\wedge dy^1_i,~ ~(\omega_2)_p =
\sum_{i=1}^m dx_i\wedge dy^2_i,~ ~(\omega_D)_p = \sum_{i=1}^m
dy^1_i\wedge dy^2_i + \sum_{i=1}^cdz_i\wedge dw_i\] Any orthogonal
basis of $T_pX$ dual to a basis  of $1$- forms  as above  is said to
be {\em adapted} to the structure, or {\em standard}. The number $m$
is the {\em rank} of the structure.
\end{dfn}
For a more intrinsic definition of WSD manifolds the reader should
refer to  ~\cite{G2}.  Here we have chosen the quickest way to
introduce them. \\
When the forms $\omega_1,\omega_2,\omega_D$ are covariant constant with respect to the Levi-Civita connection, we speak of $2$-\ka manifolds. An example of these comes from mirror symmetry for abelian varieties.
\begin{rmk}
The form $\omega_D$ is symplectic once restricted to $\omega_1^0 +
\omega_2^0$. We have therefore that $\omega_D^{dim(X) - m}\not=0$.
\end{rmk}
\begin{dfn}
1) A WSD manifold is {\em nondegenerate} if $dim(\omega_1^0\cap\omega_2^0)_p = 0$ at all points (equivalently if its dimension is $3$ times the rank).\\
2) A  WSD manifold is {\em self-dual} (SD manifold for brevity) if
all the leaves of the distribution $\omega_1^0+\omega_2^0$ have
volume one (with respect to the volume form induced by the metric)
\end{dfn}
Using Self dual manifolds, you can give a first na\"ive geometric definition of Mirror Symmetry as follows:\\
\textit{Two Calabi-Yau manifolds with $B$-field $(M,B_{M})$ and
$(\hat{M},B_{\hat{M}})$ are mirror dual if there is a Self-dual
manifold $X$ together with surjections $\pi:X \to M$ and
$\hat{\pi}:X\to \hat{M}$ such that:\\
a) $\pi^*(\omega_M) = \omega_1$,~ $\hat{\pi}^*(\omega_{\hat{M}}) = \omega_2$.\\
b) The leaves of $\omega_1^\perp$ are the fibres of $\hat{\pi}$\\
c) The leaves of $\omega_2^\perp$ are the fibres of $\pi$\\
d) The induced $B$-fields on $M$ and $\hat{M}$ are the ones
given.}\\
Here make their first appearence the B-fields $B_M$ and
$B_{\hat{M}}$, which are flat unitary gerbes on $M$ and $\hat{M}$
respectively, and which are not relevant for the discussions of this
paper. In \cite{G2} it was shown that this picture works well in the
case of elliptic curves, and for some other flat situations.\\

\textbf{Physical motivation.} One of the reasons to introduce SD manifolds however was to get rid
of special lagrangian fibrations, which are so difficult to
construct, and to be able to attack the problem of Mirror Symmetry
also when these fibrations are not expected to exist. In this more
general context one expects that the Mirror Symmetry phaenomenon
will not be obtained directly from fibrations of a SD manifold to
the dual Calabi-Yaus, but via a more sophisticated procedure, which
involves a Gromov-Hausdorff type of limit. In \cite{G3} it was shown
that for the family of anticanonical divisors in complex projective
space one can build a (real) two-dimensional family of WSD
manifolds, which degenerate in a normalized Gromov-Hausdorff sense
to the correct limits of the mirror dual Calabi-Yaus.  The picture
is the following:

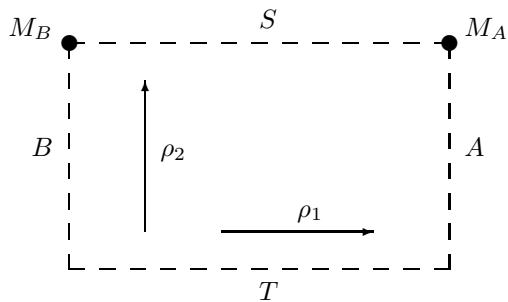
\begin{center}
\setlength{\unitlength}{1 cm}
\begin{picture}(7,5)
\put(1,1){\dashbox{0,2}(5,3)}
\put(0.2,4.1){$M_B$} \put(6.2,4.1){$M_A$} \put(3.5,4.2){$S$}
\put(3.5,0.6){$T$} \put(0.5,2.5){$B$} \put(6.2,2.5){$A$}
\put(1,4){\circle*{0,2}} \put(6,4){\circle*{0,2}}
\put(2,1.5){\vector(0,1){2}} \put(2.2,2.5){$\rho_2$}
\put(3,1.5){\vector(1,0){2}} \put(4,1.7){$\rho_1$}
\end{picture}
\end{center}
where $M_A$ and $M_B$ are the large K\"ahler and large complex
structure limits of $M$ and $\hat{M}$ respectively. To be precise,
the manifolds which come out of the costruction  of \cite{G3} are 11
dimensional (degenerate) Weakly self-dual manifolds or rank $3$.
Dimension 11 is very appealing in this context from a physical point
of view, and it brings us to the motivation for the present work.

The point of view of \cite{G3} is very different from the current
one in the main literature on mathematical Mirror Symmetry: instead
of considering the fibre product $M\times_B\hat{M}$ (when it exists)
as a device for proving Mirror Symmetry for Calabi-Yaus, the
limiting Calabi-Yaus of Mirror Symmetry are seen as very special
limits of a family of Self-Dual manifolds, which are the main
objects of study. This is actually more in line with what can be
found in the physical literature, where the $\sigma$-models defining
the string theories from which Mirror Symmetry originates are seen
as just "phases" of a unique theory, which is not necessarily in the
form of a $\sigma$-model but could very likely be similar to a
quantized Gauge theory on an 11-dimensional manifold. To make this
circle of ideas more concrete (and hence more verifiable) at the end
of \cite{G3} it is suggested that one should try to build a natural
gauge theory on Self-dual manifolds: the hope is that once quantized
this gauge theory might interpolate between the $\sigma$-models
associated to the Calabi-Yau's, and as a byproduct prove Mirror
Symmetry for them. Of course one can always put a gauge bundle on
the Self-dual manifolds "artificially", but a natural bundle which
depends only on the geometric structure would be much more
appealing. We ignore here the issue of which action to put on the
theory, but it too should be a natural geometric one. Finally, on \cite{GG}
we analyzed the situation for rank three WSD manifolds, and we found that in 
this case the corresponding natural bundle is formed by complex Lie superalgebras. 
We were able to find a geometrically motivated real form, and to split it into
simple factors. The results of \cite{GG} confirm the suspicion that on a WSD
manifold of high enough rank there could be enough natural algebraic bundles
of operators to build interesting gauge theories.\\

\textbf{The construction of $\gg$.}
From a physical point of view the case of Calabi-Yau threefolds
(i.e. rank three WSD manifolds) or fourfolds (i.e. rank four WSD manifolds) would be the most interesting one to
start with. However, its technical difficulty convinced us to start
more modestly from the case of Calabi-Yau two-folds (i.e. K3
surfaces) which correspond to rank two Self-dual manifolds. We also
considered only orientable nondegenerate Self-dual manifolds of rank
two, hence of dimension $6$. This could be considered a proof of
concept from a physicist's point of view, however Mirror Symmetry
for K3's is in itself very interesting mathematically, so we hope
that our results could have some useful geometric consequences. The rank three case is treated in our subsequent \cite{GG}, as mentioned in the previous section of this introduction. The
main result of
the present paper is the following (which is a geometric restatement of  Theorem \ref{teo:kzero}): \\

\textit{The Lie algebra $\slsc$ acts via canonical operators
(depending only
 on the geometric structure) on the smooth differential forms of any orientable
 nondegenerate WSD manifold of rank $2$.}\\

 This action generalizes naturally the action of $\sl2c$ on smooth differential
 forms of any almost \ka manifold, and is induced by a bundle action on
  the exterior power of the cotangent bundle. \\
  
Recall that a Weakly self-dual manifold is a Riemannian manifold
with three "compatible" closed differential forms. We will build a
Lie algebra of pointwise operators on complex differential forms on
$X$, as smooth sections of a bundle of Lie algebras of operators on
the complexified cotangent bundle of $X$. To start, one can define
the following operators:
\begin{dfn}
\label{dfn:Loperators} For $\phi\in \Omega^*_\C X$,
\[L_0(\phi) = \omega_D\wedge \phi,\qquad L_1(\phi) = -\omega_2\wedge
\phi,\qquad L_2(\phi) = \omega_1\wedge \phi\]
\end{dfn}
One can notice immediately the strong resemblance of the operators
above with the Lefschetz operator of K\"ahler geometry. Indeed, one
can elaborate on this similarity, and use the metric to define the
adjoints $\Lambda_j = L_j^*$ (using a pointwise procedure, as in the almost \ka case).\\
Simply using the $L_j$ and the $\Lambda_j$, one can show that the
algebra generated is isomorphic to $SL(4,\C)$ (\cite{G2}). However,
there are other natural differential forms on a WSD manifold  (which do not have a counterpart in the \ka case), namely
the volume forms of the distributions $\omega_1^\perp$,
$\omega_2^\perp$, $\omega_D^\perp$ of vectors which contract to zero
with the forms $\omega_1,\omega_2$ and $\omega_D$ respectively. If
one calls $V_0,V_1,V_2$ the corresponding wedge operators, and
$A_0,A_1,A_2$ their adjoints, the complexity of the calculations to
describe the generated Lie algebra grows a lot. 
We called $\ggr$ the algebra generated by the $L_j,V_j$ and their adjoints, and $\gg$ its complexification. To study $\gg$ 
we introduced an operator $J$, which is a complex structure on each
of the two-dimensional distributions mentioned above and generates a group isomorphic to $\mathbf{SO}(2,\R)$ (recall that we
are in the "hyperkahler" case, corresponding to Mirror Symmetry for
K3's, so an "extra" complex structure shouldn't be surprising; moreover the holonomy of a WSD  manifold in which all $\omega_1,\omega_2,\omega_D$ are invariant  is actually always included in the group generated by $J$). One
checks that all the operators introduced commute with it:
\[\forall j ~[L_j,J] = [\Lambda_j,J] = [V_j,J] = [A_j,J] = 0\]
and therefore one can try to decompose 
$\Lambda^*T^*_\C X$ with respect
to $J$ and then use Shur's Lemma to reduce to the study of the
operators on the isotypical components. One should mention that in
the (very) good cases (for instance $2$-\ka manifolds)
the operators above are all covariant constant with respect
to the metric connection, and define an action on the cohomology of
$X$ much in the same way as in the K\"ahler setting the operators
$L$ and $\Lambda$ do (due to Hodge-type identities). We don't
explore this aspect here, although it may be relevant to the
(homological) mirror map construction.
\\
Coming back to the construction,  we point out the inclusion of the Lie algebra $\gg$  inside a copy of the Clifford algebra $\mathbf{Cl}_{6,6}$. \\
Using this Clifford algebra one can identify "degree two" or "quadratic" operators (in a way similar to the ones involved in the Spinor representations on standard Spin manifolds) and among these the $\mathbf{SO}(2,\R)$-invariant ones.  A posteriori, it turns out that the operators of $\gg\oplus <J>$ за are all the $J$-invariant operators of "degree two", and this strengthens the rationale in our selection of natural operators.\\
As a last step one finds that inside
$\Lambda^*T^*X$ there is an $\mathbf{SO}(2,\R)$-isotypical component of dimension $6$,
and by direct computation we prove that  indeed the
operators restricted to this sub-representation determine a copy of
$\slsc$ (with the defining representation). 
Using the bound on the dimension of $\gg$ obtained computing "quadratic" invariants, one then shows that the representation on this isotypical component is faithful. 
This provides as a byproduct a method for giving presentation of standard Serre generators of $\gg$, explicitely written in terms of the natural geometrical
generators.

\section{Basic operators}
\label{sec:basicoperators} In this section we fix a point $p$ in the
WSD manifold $X$. The WSD structure splits the cotangent space as
$T_p^*X = W_0 \oplus W_1 \oplus W_2$ where the $W_j$ are  three
mutually orthogonal canonical distributions defined as:
\[W_0 = \{\phi\in T_p^*X~|~\phi\wedge \omega_1^2 = \phi\wedge
\omega_2^2 = 0\}\]
\[W_1 = \{\phi\in T_p^*X~|~\phi\wedge \omega_1^2 = \phi\wedge
\omega_D^2 = 0\}\]
\[W_2 = \{\phi\in T_p^*X~|~\phi\wedge \omega_2^2 = \phi\wedge
\omega_D^2 = 0\}\]

The WSD structure also determines canonical pairwise linear
identifications among $W_0,W_1$ and $W_2$, so that one can also
write $T_p^*X = W_0\otimes_\R \R^3$ or more simply
\[T_p^*X = W\otimes_\R \R^\dime\] where $W = W_0 \cong W_1 \cong
W_2$.\\
Let us now come back to the canonical  operators $L_j$ mentioned in
the introduction: \\
\textbf{Definition~\ref{dfn:Loperators}}\textit{ For $\phi\in
\Omega^*_\C X$,
\[L_0(\phi) = \omega_D\wedge \phi,\qquad L_1(\phi) = -\omega_2\wedge
\phi,\qquad L_2(\phi) = \omega_1\wedge \phi\] }\\
We now choose a (non-canonical) orthonormal basis
$\gamma_1,\gamma_2$ for $W_0$, and this together with the standard
identifications of the $W_j$ determines an orthonormal basis for
$T_p^*X$, which we write as $\{v_{ij}= \gamma_i\otimes
e_j~|~i=1,2,~j=0,1,2\}$.  We remark that the $v_{ij}$ are an
\textit{adapted coframe} for the WSD structure, and therefore we
have the explicit expressions:
\[\omega_1 = v_{10}\wedge v_{11} + v_{20}\wedge v_{21}\]
\[\omega_2 = v_{10}\wedge v_{12} + v_{20}\wedge v_{22}\]
\[\omega_D = v_{11}\wedge v_{12} + v_{21}\wedge v_{22}\]
A different choice of the $\gamma_1,\gamma_2$ would be related to
the previous one by an element in $\mathbf{O}(2,\R)$ or, taking into
account the orientability of $X$ mentioned in the Introduction, an
element of $\mathbf{SO}(2,\R)$. The Lie algebra of the group
$\mathbf{SO}(2,\R)$ expressing the change from one oriented adapted
basis to another is generated (point by point) by the global
operator $J$:
\begin{dfn}
\label{dfn:so2R} The operator $J\in End_\R(\Omega^*(X))$ is induced
by its pointwise action on the $\Lambda^*T^*_pX$ for varying $p\in
X$, defined in terms of the standard basis $v_{ij}$ as
\[J(v_{1j}) = v_{2j},\qquad J(v_{2j}) = - v_{1j}\qquad \text{for}~
j \in \{0,1,2\}\] and $J(v\wedge w) = J(v)\wedge w + v \wedge J(w)$
for $v,w \in \Lambda^*T^*_pX$
\end{dfn}
\begin{rmk}
As $J$ commutes with itself, it is well defined, independently of
the choice of an oriented adapted basis.
\end{rmk}

Using the chosen (orthonormal) basis, one can define corresponding (non canonical) wedge and contraction operators:
\begin{dfn} Let $i\in\{1,2\}$ and $j\in\{0,1,2\}$. The operators $E_{ij}$ and $I_{ij}$ are respectively the wedge and the contraction operator with the form $v_{ij}$ on $\bigwedge^*T^*X$ (defined using the given basis); we use the notation $\frac{\partial}{\partial v_{ij}}$ to indicate the element of $T_pX$ dual to $v_{ij}\in T^*_pX$:
\[E_{ij}(\phi) = v_{ij}\wedge \phi,\qquad I_{ij}(\phi) = \frac{\partial}{\partial v_{ij}}\contr
\phi\]
\end{dfn}
\begin{pro}
\label{pro:clifrelations}
The operators $E_{ij},I_{ij}$ satisfy the following relations:
\[ \forall i,j,k,l\qquad E_{ij}E_{kl} = -E_{kl}E_{ij},\quad I_{ij}I_{kl} = - I_{kl}I_{ij}\]
\[\forall i,j \qquad E_{ij}I_{ij} + I_{ij}E_{ij} = Id\]
\[\forall (i,j)\not= (k,l) \qquad E_{ij}I_{kl} = - I_{kl}E_{ij}\]
\[\forall i,j\qquad E_{ij}^* = I_{ij},\quad I_{ij}^* = E_{ij}\]
where $*$ is adjunction with respect to the metric.
\end{pro}
\dimo The proof is a simple direct verification, which we omit. \qed\\
It is then immediate to verify that:
\begin{pro}
\label{pro:dfnJ} $J$  can be expressed as
\[J = \sum_{j=0}^2\left(E_{2j}I_{1j} - E_{1j}I_{2j}\right)\]
on the whole $\bigwedge^*T^*_pX$. From this expression and the previous proposition one obtains that  $J^* = -J$, i.e.  for every $p$ the Lie algebra generated
by $J$ is a subalgebra of $\mathbf{o}(\bigwedge^*T^*_pX)$ isomorphic to $\so2\cong \R$. Moreover, the exponential images inside $Aut_\R(\Omega^*(X))$of the operators of type $tJ$ for $t\in \R$ form a group isomorphic to $\mathbf{SO}(2,\R) \cong S^1$, as this isomorphism holds for the (faithful) restriction of the group action to $T^*_pX$.
\end{pro}
Using the (non canonical) operators $E_{ij}$ we can obtain simple expressions for the
pointwise action of the other canonical operators, the volume forms $V_j$:
\begin{dfn}
\label{dfn:Voperators} For $\phi\in \bigwedge^*T^*_pX$,
\[V_0(\phi) = E_{10} E_{20}(\phi),\qquad
V_1(\phi) = E_{11} E_{21}(\phi),\qquad V_2(\phi) =
E_{12} E_{22}(\phi)\]
\end{dfn}
Remember however that the operators $V_j$ do not depend on the choice of a basis, as they
are simply multiplication by the volume forms of the spaces $W_j$.

We use the $v_{ij}$ also as a orthonormal basis for the complexified
space $T_p^*\otimes_\R\C$ (with respect to the induced hermitian
inner product). We indicate with the same symbols $V_j$ the
complexified operators acting on the spaces $\bigwedge_\C^*T^*_pX$.

The riemannian metric induces a Riemannian metric on $T^*_p X$ and
on the space $\bigwedge^*T^*_pX$.
\begin{dfn}
\label{dfn:lambda-A-operators} For $j\in\{0,1,2\}$
\[\Lambda_j = L_j^*,\qquad A_j =
V_j^*\]
\end{dfn}
By
construction the canonical operators $L_j,V_j,\Lambda_j,A_j$ on
$\bigwedge^*T^*_pX$ are the pointwise restrictions of
corresponding global operators on smooth differential forms, which
we indicate with the same symbols: for $j \in \{0,1,2\}$,
\[L_j,V_j,\Lambda_j,A_j: \Omega^*(X)\to
\Omega^*(X)\]
Summing up:
\begin{dfn}
\label{dfn:algebragg} The $*$-Lie algebra $\ggr$ is the $*$-Lie
subalgebra of $End_\R\left(\Omega^*(X)\right)$ generated by the
operators
\[\{L_j,V_j,\Lambda_j,A_j~|~\text{for}~ ~ j = 0,1,2\}\]
The $*$ operator on $\ggr$ is induced by the adjoint with respect to
the Riemannian metric. The $*$-Lie algebra $\gg$ is $\ggr\otimes \C$, and is in a natural way a $*$-Lie subalgebra of $End_\C\left(\Omega_\C^*(X)\right)$. The $*$ operator on $\gg$ is induced by the adjoint with respect to
the induced Hermitian metric.
\end{dfn}
 The canonical splitting
$T_p^*X = W_0 \oplus W_1 \oplus W_2$ together with the canonical
identifications $W_0\cong W_1\cong W_2$ induce an action of the
symmetric group $\mathcal{S}_3$, which propagates to
$\bigwedge^*T^*X$ and to its $\mathcal{C}^\infty$ sections. At
every point, the action can be written explicitly in terms of the
basis as
\[\sigma(v_{ij}) = v_{i\sigma(j)}\]
The induced action on endomorphisms via
conjugation, $\sigma(\phi) = \sigma\circ\phi\circ\sigma^{-1}$,
preserves $\gg$. Indeed, one can check directly using the basis
$v_{ij}$ at every point that for $\sigma \in \s3$
\[\sigma(V_j) = V_{\sigma(j)},\qquad \sigma(L_j) = \epsilon(\sigma)L_{\sigma(j)}\]
Since $\s3$ acts on $\gg$ by conjugation with unitary operators, its
action commutes with adjunction (the $*$ operator), and therefore
\[\sigma(A_j) = A_{\sigma(j)},\qquad \sigma(\Lambda_j) = \epsilon(\sigma)\Lambda_{\sigma(j)}\]
Moreover, one also has that $\sigma(J) = J$ which means that the
action of $\mathcal{S}_3$ commutes with that of $\so2$.\\
\section{The action of $\so2$}
\label{sec:so2R} When one deals with mirror simmetry for $2$-\ka manifolds (see the Introduction), the WSD manifolds which arise have the property that the
forms $\omega_1,\omega_2$ and $\omega_D$ are covariant constant with
respect to the metric. In this case, the maximal possible holonomy
of the WSD manifold $X$ is included in the $\so2$ generated by the
operator $J$. We will show now that $J$ commutes with $\gg$. Our
proof will be strictly algebraic, so that the commutativity between
$\so2$ and $\gg$ will hold also on WSD manifolds for which the
holonomy is more general.
\begin{dfn} Given $n\in\Z$, we indicate with $V_n$ the one
dimensional complex representation of $\mathbf{SO}(2,\R) \cong S^1\
\cong \R / \Z$ given by the character:
\[\theta\to e^{2\pi \imath n\theta}\]
\end{dfn}
\begin{pro}
\label{table} Under the $\mathbf{SO}(2,\R)$ representation induced by the
operator
$J$, for any $p\in X$ :\\
1) The space $\bigwedge^1(T^*_\C X_p)$ splits as \[ V_{-1}^{\oplus
3}\bigoplus V_{1}^{\oplus 3} \]\\
2) The whole space $\bigwedge^*(T^*_\C X_p)$ splits according to the
following picture:
\[
\begin{array}{ccccccccccccc}
   \bigwedge^0 (T^*_\C X_p) & = &   &   &   & V_{0} &  &   &    \\
     &  &   &   &   &  &  &   &     \\
   \bigwedge^1 (T^*_\C X_p) & = &   &   & V_{-1}^{\oplus
3} & \bigoplus  & V_{1}^{\oplus
3} &  &    \\
    &  &   &   &   &  &  &   &    \\
    \bigwedge^2 (T^*_\C X_p) & = &  & V_{-2}^{\oplus  3} & \bigoplus  &
    V_{0}^{\oplus 9} & \bigoplus & V_{2}^{\oplus 3} &   \\
        &  &   &   &   &  &  &   &    \\
    \bigwedge^3 (T^*_\C X_p) & = & V_{-3} & \bigoplus & V_{-1}^{\oplus
    9} & \bigoplus & V_{1}^{\oplus 9} & \bigoplus & V_{3}  \\
        &  &   &   &   &  &  &   &    \\
    \bigwedge^4 (T^*_\C X_p) & = &  & V_{-2}^{\oplus  3} & \bigoplus  &
    V_{0}^{\oplus 9} & \bigoplus & V_{2}^{\oplus 3} &     \\
        &  &   &   &   &  &  &   &    \\
    \bigwedge^5 (T^*_\C X_p) & =&   &   & V_{-1}^{\oplus
3} & \bigoplus  & V_{1}^{\oplus
3} &  &   \\
    &  &   &   &   &  &  &   &     \\
   \bigwedge^6 (T^*_\C X_p) & =&    &   &   & V_{0} &  &   &
\end{array}
\]
\end{pro}
\[ \]
\dimo 1) The space $T^*_\C X_p$ is a direct sum of the three $W_j$,
and each one of these is the standard two dimensional real
representation of $\so2$. We therefore diagonalize the
representation introducing a new basis for each $W_j =
<v_{1j},v_{2j}>$:
\[w_j = v_{1j} + \imath\ v_{2j},\qquad\overline{w}_j = v_{1j} - \imath
v_{2j}\] From the definition of $J$, one has then for every $j\in
\{0,1,2\}$
\[J(w_j) = -\imath w_j,\qquad
J(\overline{w}_j) = \imath \overline{w}_j\] Therefore one has for
every $j\in \{0,1,2\}$
\[<w_j> \cong V_{-1},\qquad <\overline{w}_j> \cong V_1\]
2) To prove the general case, we use the fact that the operator $J$
determines an almost complex structure on the manifold $X$,
compatible with the metric. From this, following standard arguments,
the complex differential forms and also the elements of
$\bigwedge^*T^*_\C X_y$ for any $y\in Y$ can be divided according to
their type:
\[\bigwedge^*T^*_\C X_y = \bigoplus_{n=0}^{dim X}\bigoplus_{p+q =
n}\bigwedge^{p,q}T^*_\C X_y\] In the notation adopted in the proof
of the first statement, one has
\[\bigwedge^{p,q}T^*_\C X_y = <w_{i_1}\wedge\cdots\wedge w_{i_p}\wedge
\overline{w}_{j_1}\wedge\cdots\wedge
\overline{w}_{j_q}~|~i_1,...,j_q~\in\{0,1,2\}>\] From the definition
of the action of $J$ one has therefore that for any $p,q$
\[\bigwedge^{p,q}T^*_\C X_y \cong V_{q-p}^{\oplus k}\]
with $k = \begin{pmatrix}3\\ p\end{pmatrix}\begin{pmatrix}3\\
q\end{pmatrix}$ from which the second statement of the proposition
can be esily deduced.
\qed\\
\begin{teo}
\label{olonomy} The operators $L_j,V_j$ for $j\in\{0,1,2\}$ commute
with the generator $J$ of $\so2$.
\end{teo}
\dimo We prove the statements by a direct computation using the
basis $v_{ij}$; moreover, using the action of $\mathcal{S}_3$ (which
permutes the $L_j,V_j$ and fixes $J$), it is enough to prove the
commutativity for $L_0$ and $V_0$. It useful to rewrite $\omega_0$
(and hence $L_0$ which is wedge with $\omega_0$) in terms of the
basis generated by the $w_j$:
\[\omega_0 = v_{11}\wedge v_{12} + v_{21}\wedge v_{22} =
\frac{1}{2}\left(w_{1}\wedge \overline{w}_{2} - w_{2}\wedge
\overline{w}_{1}\right)\] and then:
\[[J,L_0](w_{i_1}\wedge\cdots\wedge w_{i_p}\wedge
\overline{w}_{j_1}\wedge\cdots\wedge \overline{w}_{j_q}) =\]
\[J\left(\frac{1}{2}\left(w_{1}\wedge \overline{w}_{2} - w_{2}\wedge
\overline{w}_{1}\right)\right)\wedge
\left(w_{i_1}\wedge\cdots\wedge w_{i_p}\wedge
\overline{w}_{j_1}\wedge\cdots\wedge \overline{w}_{j_q}\right) + \]
\[\left(\frac{1}{2}\left(w_{1}\wedge \overline{w}_{2} - w_{2}\wedge
\overline{w}_{1}\right)\right)\wedge
J\left(w_{i_1}\wedge\cdots\wedge w_{i_p}\wedge
\overline{w}_{j_1}\wedge\cdots\wedge \overline{w}_{j_q}\right) - \]

\[- \frac{1}{2}\left(w_{1}\wedge \overline{w}_{2} - w_{2}\wedge
\overline{w}_{1}\right)\wedge J(w_{i_1}\wedge\cdots\wedge
w_{i_p}\wedge \overline{w}_{j_1}\wedge\cdots\wedge
\overline{w}_{j_q})\] Therefore the result follows from the fact
that
\[J(\frac{1}{2}\left(w_{1}\wedge \overline{w}_{2} - w_{2}\wedge
\overline{w}_{1}\right) = 0\] 
as $w_j$ and  $\overline{w}_k$ have opposite weight with respect to $J$ for any $j,k$.\\
Similarly,
$[J,V_0] = 0$ follows from the fact that for any $\alpha$
\[V_0(\alpha) = v_{10}\wedge v_{20}\wedge \alpha = \frac{\imath}{2}\ w_{0}\wedge \overline{w}_{0}\wedge
\alpha\] \qed \\
From the previous theorem  one obtains the
following corollary, which holds on any WSD manifold (not necessarily $2$-\ka):
\begin{cor} The algebra $\gg$ commutes with the action of $\so2$
induced by $J$.
\end{cor}
\dimo We already know that $[J,L_j] = [J,V_j] = 0$ for $j\in
\{0,1,2\}$. The corresponding commutation relations for the adjoint
generators   $\Lambda_j,A_j$ of $\gg$ follow from the fact that $J^*
= -J$, as noticed in Proposition \ref{pro:dfnJ}. \qed \\
\begin{rmk}
\label{rmktable} From Schurs's lemma it follows that the columns of
the diagram of Proposition \ref{table} are preserved by the action
of $\gg$.
\end{rmk}
\section{An irreducible representation of $\gg$}
\label{sec:irreducible} Looking at the table in Proposition \ref{table}  we notice that the second column from the left is a representation of $\gg$ (by Remark ~\ref{rmktable}) of dimension $6$:
\[V \cong V_{-2}^{\oplus
6} = <w_{0}\wedge w_{1},~w_{0}\wedge w_{2},~ w_{1}\wedge w_{2},
~w_{0}\wedge w_{1} \wedge w_{2}\wedge \overline{w}_{0},\]
\[w_{0}\wedge
w_{1} \wedge w_{2}\wedge \overline{w}_{1},~ w_{0}\wedge w_{1} \wedge
w_{2}\wedge \overline{w}_{2}>\] In this section we will compute
explicitely this representation.

Using the above described basis, it is not difficult to
compute the matrices by hand:
\begin{pro}
\label{teo:representation} Indicating with $\beta$ the ordered basis for $V$
indicated above, the matrices for the  (restrictions to $V$ of) the
generators of $\gg$ are the following:
\[\small{M_\beta(L_0) ~= ~\begin{pmatrix}
0 & 0 & 0 & 0 & 0 & 0\\
0 & 0 & 0 & 0 & 0 & 0\\
0 & 0 & 0 & 0 & 0 & 0\\
0 & 0 & 0 & 0 & 0 & 0\\
-\frac{1}{2} & 0 & 0 & 0 & 0 & 0\\
0 & -\frac{1}{2} & 0 & 0 & 0 & 0\end{pmatrix},
\quad M_\beta(\Lambda_0) ~= ~\begin{pmatrix}
0 & 0 & 0 & 0 & -2 & 0\\
0 & 0 & 0 & 0 & 0 & -2\\
0 & 0 & 0 & 0 & 0 & 0\\
0 & 0 & 0 & 0 & 0 & 0\\
0 & 0 & 0 & 0 & 0 & 0\\
0 & 0 & 0 & 0 & 0 & 0\end{pmatrix}}\]
\[\small{M_\beta(L_1) ~= ~\begin{pmatrix}
0 & 0 & 0 & 0 & 0 & 0\\
0 & 0 & 0 & 0 & 0 & 0\\
0 & 0 & 0 & 0 & 0 & 0\\
\frac{1}{2} & 0 & 0 & 0 & 0 & 0\\
0 & 0 & 0 & 0 & 0 & 0\\
0 & 0 & -\frac{1}{2} & 0 & 0 & 0\end{pmatrix}, \quad
M_\beta(\Lambda_1) ~= ~\begin{pmatrix}
0 & 0 & 0 & 2 & 0 & 0\\
0 & 0 & 0 & 0 & 0 & 0\\
0 & 0 & 0 & 0 & 0 & -2\\
0 & 0 & 0 & 0 & 0 & 0\\
0 & 0 & 0 & 0 & 0 & 0\\
0 & 0 & 0 & 0 & 0 & 0\end{pmatrix}}\]
\[\small{M_\beta(L_2) ~= ~\begin{pmatrix}
0 & 0 & 0 & 0 & 0 & 0\\
0 & 0 & 0 & 0 & 0 & 0\\
0 & 0 & 0 & 0 & 0 & 0\\
0 & \frac{1}{2} & 0 & 0 & 0 & 0\\
0 & 0 & \frac{1}{2} & 0 & 0 & 0\\
0 & 0 & 0 & 0 & 0 & 0\end{pmatrix}, \quad M_\beta(\Lambda_2) ~=
~\begin{pmatrix}
0 & 0 & 0 & 0 & 0 & 0\\
0 & 0 & 0 & 2 & 0 & 0\\
0 & 0 & 0 & 0 & 2 & 0\\
0 & 0 & 0 & 0 & 0 & 0\\
0 & 0 & 0 & 0 & 0 & 0\\
0 & 0 & 0 & 0 & 0 & 0\end{pmatrix}}\]
\[\small{M_\beta(V_0) ~= ~\begin{pmatrix}
0 & 0 & 0 & 0 & 0 & 0\\
0 & 0 & 0 & 0 & 0 & 0\\
0 & 0 & 0 & 0 & 0 & 0\\
0 & 0 & \frac{\imath}{2} & 0 & 0 & 0\\
0 & 0 & 0 & 0 & 0 & 0\\
0 & 0 & 0 & 0 & 0 & 0\end{pmatrix}, \quad M_\beta(A_0) ~=
~\begin{pmatrix}
0 & 0 & 0 & 0 & 0 & 0\\
0 & 0 & 0 & 0 & 0 & 0\\
0 & 0 & 0 & -2\imath & 0 & 0\\
0 & 0 & 0 & 0 & 0 & 0\\
0 & 0 & 0 & 0 & 0 & 0\\
0 & 0 & 0 & 0 & 0 & 0\end{pmatrix}}\]
\[\small{M_\beta(V_1) ~= ~\begin{pmatrix}
0 & 0 & 0 & 0 & 0 & 0\\
0 & 0 & 0 & 0 & 0 & 0\\
0 & 0 & 0 & 0 & 0 & 0\\
0 & 0 & 0 & 0 & 0 & 0\\
0 & -\frac{\imath}{2} & 0 & 0 & 0 & 0\\
0 & 0 & 0 & 0 & 0 & 0\end{pmatrix}, \quad M_\beta(A_1) ~=
~\begin{pmatrix}
0 & 0 & 0 & 0 & 0 & 0\\
0 & 0 & 0 & 0 & 2\imath & 0\\
0 & 0 & 0 & 0 & 0 & 0\\
0 & 0 & 0 & 0 & 0 & 0\\
0 & 0 & 0 & 0 & 0 & 0\\
0 & 0 & 0 & 0 & 0 & 0\end{pmatrix}}\]
\[\small{M_\beta(V_2) ~= ~\begin{pmatrix}
0 & 0 & 0 & 0 & 0 & 0\\
0 & 0 & 0 & 0 & 0 & 0\\
0 & 0 & 0 & 0 & 0 & 0\\
0 & 0 & 0 & 0 & 0 & 0\\
0 & 0 & 0 & 0 & 0 & 0\\
\frac{\imath}{2} & 0 & 0 & 0 & 0 & 0\end{pmatrix}, \quad
M_\beta(A_2) ~= ~\begin{pmatrix}
0 & 0 & 0 & 0 & 0 & -2\imath\\
0 & 0 & 0 & 0 & 0 & 0\\
0 & 0 & 0 & 0 & 0 & 0\\
0 & 0 & 0 & 0 & 0 & 0\\
0 & 0 & 0 & 0 & 0 & 0\\
0 & 0 & 0 & 0 & 0 & 0\end{pmatrix}}\]

\end{pro}
\dimo Direct computation using the  basis generated by the $w_j$.
\qed
\begin{cor}
\label{cor:sl6} The algebra generated by the restriction of $\gg$ to
$V$ is isomorphic to $\slsc$, with $V$ its natural representation.
\end{cor}
One can sum up the computations above in the following theorem:
\begin{teo}
\label{teo:restriction} There is an exact sequence of Lie algebras
\[0 \to K\to \gg \to \slsc\to 0\]
given by the restriction to $V$. 
\end{teo}
In the next section we will prove that $K = \{0\}$, and therefore the representation $V$ is faithful and $\gg \cong \slsc$.
\section{Quadratic invariants}
\label{sec:quadratic}
We begin by showing that the action of Lie algebra $\gg$  is induced by a (non-canonical) Clifford algebra representation. We use for simplicity the canonical identification $T^{**}X_p \cong TX_p$ without further comment, so that if $\{v_{ij}\}$ is a basis for $T^*_pX$, then $\{\frac{\partial}{\partial v_{ij}}\}$ is the corresponding dual basis for $T_pX$.
\begin{dfn}
For $p\in X$, the Clifford algebra $\mathcal{C}_p$ is 
\[\mathcal{C}_p = Cl(T_pX \oplus T^*_pX,q)\]
with the quadratic form $q$ induced by the  metric
\[\begin{array}{ll}\forall i,j,h,k & <v_{ij},v_{hk}> = 0\\
\forall i,j,h,k & <\frac{\partial}{\partial v_{ij}},\frac{\partial}{\partial v_{hk}}> = 0\\
\forall (i,j)\not= (h,k) &  <v_{ij},\frac{\partial}{\partial v_{hk}}> = 0\\
\forall i,j & <v_{ij},\frac{\partial}{\partial v_{ij}}> = -\frac{1}{2}
\end{array}\]
\end{dfn}
\begin{rmk} The Clifford algebras $\mathcal{C}_p$ for varying $p$ define a Clifford bundle $\mathcal{C}$ on $X$, as the definition of $\mathcal{C}_p$ is independent on the choice of a basis. Indeed, the quadratic form used to define it is simply induced by $-\frac{1}{2}$ times the natural bilinear pairing  $T_pX \otimes T^*_pX\to \R$.
\end{rmk}
\begin{pro}
The Clifford algebra $\mathcal{C}_p$ has a canonical representation $\rho_p$ on $\bigwedge T^*_p X$, induced by the operators $E_{ij}$ and $I_{ij}$ via the map
\[\rho_p(v_{ij}) =  E_{ij},\qquad \rho_p\left(\frac{\partial}{\partial v_{ij}}\right) =  I_{ij}\]
\end{pro}
\dimo The Clifford relations
\[\phi\psi + \psi\phi =  -2<\phi,\psi>\]
are precisely the content of Proposition \ref{pro:clifrelations}. The representation is canonical, even if the operators $E_{ij}$ and $I_{ij}$ are not, because it can be defined in a basis independent way as
\[\rho_p(v)(\alpha) = v\wedge \alpha,\qquad \rho_p\left(\frac{\partial}{\partial v}\right) = \frac{\partial}{\partial v}\contr \alpha\]
\qed\\
Abusing slightly the notation, we will identify $\mathcal{C}_p$ with its (faithful) image inside $End_\R\left(\bigwedge^*T^*_pX\right)$, and we will omit any reference to the map $\rho_p$. Actually, as the representation above is a real analogue of the Spinor representation, it is easy to check that the map $\rho_p$  is an isomorphism of associative algebras. One then has:
\begin{dfn}
\label{dfn:quadratic} The linear subspace $\mathcal{C}^2_p$ of $\mathcal{C}_p$ is the image of the natural map $\bigwedge^2(T_pX\oplus T_p^*X)\to \mathcal{C}_p$. The linear subspace $\mathcal{C}^0_p$ of $\mathcal{C}_p$ is the subspace generated by $1$.
\end{dfn}
Recall that $\mathcal{C}^2_p$ is a Lie subalgebra  of $\mathcal{C}_p$ (with the commutator bracket).
\begin{pro} The Lie algebra $\ggr_p$ and the operator $J$ sit inside $\mathcal{C}^2_p$ for all $p \in X$.
\label{pro:subclif}
\end{pro}
\dimo The operators $L_j$, the $\Lambda_j$, the $V_j$ and the $A_j$ lie inside $\mathcal{C}^2_p\oplus \mathcal{C}^0_p$ by Proposition \ref{pro:clifrelations} and the fact that $\omega_1,\omega_,\omega_D$ lie in $\bigwedge^2T^*_pX$.  The operator $J$ lies inside $\mathcal{C}^2_p\oplus \mathcal{C}^0_p$ by Proposition \ref{pro:dfnJ}. By definition the elements $\mathcal{C}^2_p$  are commutators, and therefore have trace zero in any representation, and hence also in the $\rho_p$. Moreover, again by inspection all the generators of $\ggr_p$ have trace zero once represented via $\rho_p$ (they are nilpotent), and therefore they must lie inside $\mathcal{C}^2_p$. The operator $J$ is in the Lie algebra of the isometry group, and therefore it too has trace zero and hence sits inside $\mathcal{C}^2_p$. As $\mathcal{C}^2_p$ is closed under the commutator bracket of $\mathcal{C}_p$, and this commutator coincides with the composition bracket of operators, we have the conclusion.
\qed
\begin{rmk}
Giving degree $1$ to the operators $E_{ij}$ and degree $-1$ to the operators $I_{ij}$, we induce a $\Z$-degree on  $\mathcal{C}_p$. This degree coincides with the degree of the operators induced from the grading on  the forms from $\bigwedge^*T^*X$.
\end{rmk}
\begin{rmk} For any $p\in X$, the Clifford algebra $\mathcal{C}_p$ is isomorphic to $\mathbf{Cl}_{6,6}$, as the metric used to define it has signature $(6,6)$. The previous proposition therefore shows that $\ggr_p$ is a Lie subalgebra of $\mathbf{Cl}_{6,6}^2 \cong \mathbf{spin}(6,6) = \mathbf{so}(6,6)$, generated by smooth global sections of the Clifford bundle $\mathcal{C}$.
\end{rmk}
The operator $J$ acts on all of $\mathcal{C}_p$ by adjunction with respect to the commutator bracket, and sends its quadratic part $\mathcal{C}^2_p$  to itself from Proposition \ref{pro:subclif}. \\
We will show that the space of $J$-invariants inside $\mathcal{C}^2_p$ (the ``quadratic" $J$-invariants) coincides with $\gg$. To describe it explicitely, let us introduce the following notation:
\begin{dfn}
\[E_{w_j} = E_{1j} + \imath E_{2j},\qquad E_{\overline{w}_j} = E_{1j} - \imath E_{2j}\]
\[I_{w_j} = I_{1j} - \imath I_{2j},\qquad I_{\overline{w}_j} = I_{1j} + \imath I_{2j}\] 
\end{dfn}
\begin{lem} The adjoint action of the operator $J$ on $E_{w_j},I_{w_j},E_{\overline{w}_j},I_{\overline{w}_j}$ is:
\[[J,E_{w_j}] = -\imath E_{w_j},\quad [J,I_{w_j}] = \imath I_{w_j}\]
\[[J,E_{\overline{w}_j}] = \imath E_{\overline{w}_j},\quad[J,I_{\overline{w}_j}] = -\imath I_{\overline{w}_j}\]
\end{lem}
\dimo It is enough to consider the corresponding $J$-weights of the $w_j,\overline{w}_j$. \qed
\begin{pro}
The following $36$ operators provide a linear basis for the quadratic $J$-invariants:
\begin{enumerate} 
\item $[E_{w_0},E_{\overline{w}_1}], [E_{w_0},E_{\overline{w}_2}],[E_{w_1},E_{\overline{w}_2}],[E_{w_1},E_{\overline{w}_0}],[E_{w_2},E_{\overline{w}_0}],[E_{w_2},E_{\overline{w}_1}]$
\item $[I_{w_0},I_{\overline{w}_1}], [I_{w_0},I_{\overline{w}_2}],[I_{w_1},I_{\overline{w}_2}],[I_{w_1},I_{\overline{w}_0}],[I_{w_2},I_{\overline{w}_0}],[I_{w_2},I_{\overline{w}_1}]$ 
\item $[E_{w_0},E_{\overline{w}_0}], [E_{w_1},E_{\overline{w}_1}],[E_{w_2},E_{\overline{w}_2}]$ 
\item $[I_{w_0},I_{\overline{w}_0}], [I_{w_1},I_{\overline{w}_1}],[I_{w_2},I_{\overline{w}_2}]$ 
\item $[E_{w_0},I_{w_1}],[E_{w_0},I_{w_2}],[E_{w_1},I_{w_0}],[E_{w_1},I_{w_2}],[E_{w_2},I_{w_0}],[E_{w_2},I_{w_1}]$
\item $[E_{\overline{w}_0},I_{\overline{w}_1}],[E_{\overline{w}_0},I_{\overline{w}_2}],[E_{\overline{w}_1},I_{\overline{w}_0}],[E_{\overline{w}_1},I_{\overline{w}_2}],[E_{\overline{w}_2},I_{\overline{w}_0}],[E_{\overline{w}_2},I_{\overline{w}_1}]$ 
\item $[E_{w_0},I_{w_0}], [E_{w_1},I_{w_1}], [E_{w_2},I_{w_2}], [E_{\overline{w}_0},I_{\overline{w}_0}], [E_{\overline{w}_1},I_{\overline{w}_1}], [E_{\overline{w}_2},I_{\overline{w}_2}]$
\end{enumerate}
\end{pro}
\dimo The $J$-weight of a bracket of $J$-homogeneous operators is the sum of the respective weights. The quadratic "monomials" (with respect to the bracket) in the $E_{w_j},I_{w_j},E_{\overline{w}_j},I_{\overline{w}_j}$ are all $J$-homogeneous, and therefore to find a basis of $J$-invariant quadratic operators it is enough to identify the $J$-invariant quadratic monomials. To be $J$-invariant means simply to have weight zero, and the computation of the $J$-weight of the quadratic mononials follows immediately from those of $E_{w_j},I_{w_j},E_{\overline{w}_j},I_{\overline{w}_j}$, which are respectively $-\imath,\imath,\imath,-\imath$. \qed

We end this section with the following:
\begin{teo}
\label{teo:kzero}
In the exact sequence of Theorem \ref{teo:restriction} the kernel $K$ is equal to $\{0\}$. The algebra $\gg$ is therefore isomorphic to $\slsc$.
\end{teo}
\dimo Since $\gg$ is included in the Lie algebra of quadratic invariants, it is enough to show that $J\not\in \gg$, as from this and the previous proposition it follows that  $dim_\C(\gg)\leq 35$. As $\gg$ maps surjectively to $\slsc$ which has dimension $35$, the kernel must be zero. When restricted to the subrepresentation $V$, the generators of $\gg$ have all trace zero by inspection of their matrices. However, by definition of $V$, $J$ restricted to it is multiplication by $-2\imath$, and has therefore trace equal to $-12\imath$.\qed
\begin{cor} The Lie  algebra $\gg\oplus <J>$ equals the Lie algebra of quadratic invariants inside $\mathcal{C}^2_p$.
\end{cor}  
\section{A geometric presentation of  Serre generators}
\label{sec:geompres} 
In this section, to gain a better geometric understanding of the representation $\gg$ of $\slsc$, we explore in greater detail its relation to the geometric structure of a WSD manifold. In particular, we give a presentation of a natural choice of Cartan subalgebra and Serre generators in terms on the  
geometric generators $L_j,\Lambda_j,V_j,A_j$. 

The $L_j$ operators are similar in
nature to the Lefschetz operators of a K\"ahler manifold. This analogy is what provided the initial interest in the algebraic structure of $\gg$. Similarly to
 the corresponding standard construction of a representation of
$\mathbf{sl}(2,\C)$, we define
\begin{dfn}
\label{dfn:Hoperators}  For $j\in\{0,1,2\}$
\[H_j = [L_j,\Lambda_j]\]
\end{dfn}
These operators are self-adjoint, as $L_j^* = \Lambda_j$ by
definition. As in the context of K\"ahlerian geometry, for every $j$
the algebra $<L_j,\Lambda_j,H_j>$ turns out to be a copy of $\sl2c$.
Moreover, the following proposition shows that the operators $H_j$
are semisimple on the whole algebra $\gg$, and therefore generate a
toral subalgebra of $\gg$:
\begin{pro} \label{pro:weights}
The geometric operators $H_j$
generate a toral subalgebra of $\gg$, and the following relations hold: for $j\not= k \in \{0,1,2\}$
\begin{enumerate}
\item $[H_j,L_j] = 2L_j,\qquad [H_j,\Lambda_j] = -2\Lambda_j$
\item $[H_j,L_k] = L_k,\qquad [H_j,\Lambda_k] = -\Lambda_k$
\item $[H_j,V_j] = 0,\qquad [H_j,A_j] = 0$
\item $[H_j,V_k] = 2V_k,\qquad [H_j,A_k] = -2A_k$
\end{enumerate}
\end{pro}
\dimo In view of Theorem \ref{teo:kzero}, at this point the quickest method of proof of this proposition is to refer to the explicit matrices of the (faithful) restriction of $\gg$ to $V$.  \qed\\
The whole algebra $\gg$ splits into a direct sum of weight spaces
with respect to $<H_0,H_1,H_2>$, as this subalgebra is toral. The weight of $L_0$ with respect
to the basis dual to $H_0,H_1,H_2$ is:
\[\alpha_{L_0}~=
~(\alpha_{L_0}(H_0) , \alpha_{L_0}(H_1) , \alpha_{L_0}(H_2))~= ~( 2
,1 , 1)\] The full list is:
\[\alpha_{L_0}= (2,1,1),\qquad \alpha_{\Lambda_0}= - \alpha_{L_0}\]
\[\alpha_{L_1}= (1,2,1),\qquad \alpha_{\Lambda_1}= - \alpha_{L_1}\]
\[\alpha_{L_2}= (1,1,2),\qquad \alpha_{\Lambda_2}= - \alpha_{L_2}\]
\[\alpha_{V_0}= (0,2,2),\qquad \alpha_{A_0}= - \alpha_{V_0}\]
\[\alpha_{V_1}= (2,0,2),\qquad \alpha_{A_1}= - \alpha_{V_1}\]
\[\alpha_{V_2}= (2,2,0),\qquad \alpha_{A_2}= - \alpha_{V_2}\]
\[ \]
To find a natural geometric expression for two ad-semisimple elements which complete $<H_0,H_1,H_2>$ to a Cartan subalgebra we
look at the generators $V_j$ and $A_j$. However, it turns out that
the natural candidates $[V_j,A_j]$ already lie in the algebra
$<H_0,H_1,H_2>$. We instead build the new operators by "subtracting"
 from the $V_j$ their weight $\alpha_{V_j}$:
\begin{dfn}
\label{dfn:Soperators} We define
\[S_0 = \imath[[[V_0,\Lambda_1],\Lambda_2],L_0] \]
\[S_1 = \imath[[[V_1,\Lambda_2],\Lambda_0],L_1] \]
\[S_2 = \imath[[[V_2,\Lambda_0],\Lambda_1],L_2]\] and denote by $\mathcal{H}$ the
Lie algebra (over $\C$):
\[\h = < H_0,H_1,H_2,S_0,S_1,S_2 >\]
\end{dfn}
The coefficients $\imath$ which appear in the formulas above are dictated by the fact that  with this choice the (diagonal) matrices of the $S_j$ restricted to $V$ have integer entries.
\begin{pro}
The algebra $\h$ is a Cartan subalgebra of $\gg$. More precisely, the following are the diagonals of the operators $H_0,...,S_2$ once restricted to $V$ 
{\small
\[
H_0:\begin{pmatrix}-1\\ -1\\ 0\\ 0\\ 1\\ 1\end{pmatrix},~
H_1:\begin{pmatrix}-1\\ 0\\ -1\\ 1\\ 0\\ 1\end{pmatrix},~
H_2:\begin{pmatrix}0\\ -1\\ -1\\ 1\\ 1\\ 0\end{pmatrix},~
S_0:\begin{pmatrix}-1\\ 1\\ 0\\ 0\\ 1\\ -1\end{pmatrix},~
S_1:\begin{pmatrix}1\\ 0\\ -1\\ -1\\ 0\\ 1\end{pmatrix},~
S_2:\begin{pmatrix}0\\ -1\\ 1\\ 1\\ -1\\ 0\end{pmatrix}
\]}
\end{pro}
\dimo The computation of the matrices above shows that, once restricted to $V$, the algebra $\h$ spans the space of diagonal matrices of trace zero in the given basis. \qed
\begin{rmk}
The computation above shows also that  operators $S_0,S_1,S_2$ safisfy the relation
\[S_0 + S_1 + S_2 = 0\]
\end{rmk}
Even if from the previous proposition we know that $\h$ is maximal toral inside $\gg$, the natural geometric generators $L_j,\Lambda_j$ are not eigenvectors for the adjoint action of the $S_k$.  At this point however it is possible to single out in natural geometric terms  operators of $\gg$ which have 
"pure" weight with respect to the algebra $\h$ and which contain in their linear span the $L_j,\Lambda_j$:
\begin{dfn}
\label{dfn:lij} For $j\in\{0,1,2\}$
\[L_{1j} = -2L_j +
[S_j,L_j],\qquad L_{2j} = 2L_j + [S_j,L_j] \]
\[\Lambda_{1j} = -2\Lambda_j -
[S_j,\Lambda_j],\qquad \Lambda_{2j} = 2\Lambda_j - [S_j,\Lambda_j]
\]
\end{dfn}
\begin{pro} Indicating with $e^h_k$ the $6\times 6$ matrix with a
$1$ in position $k$ (row) and $h$ (column) and zero otherwise, the
matrices of the operators $L_{ij}$ and $\Lambda_{ij}$ restricted on
$V$ are:
\[\begin{array}{lll}
L_{1 0} = 2e_6^2 & L_{1 1} = -2e_4^1 & L_{1 2} = -2e_5^3\\
\\
L_{2 0} = -2e_5^1 & L_{2 1} = -2e_6^3 & L_{2 2} = 2e_4^2\\
\\
\Lambda_{1 0} = 8e_2^6 & \Lambda_{1 1} = -8e_1^4 & \Lambda_{1 2} = -8e_3^5\\
\\
\Lambda_{2 0} = -8e_1^5 & \Lambda_{2 1} = -8e_3^6 & \Lambda_{2 2} =
8e_2^4\end{array}\]
\end{pro}
\begin{cor}
\label{cor:serreLij}We have the following relations for the
operators of $\gg$ restricted to $V$:
\[[H_k,L_{ij}] = (1+\delta_{kj})L_{ij},
\quad [H_k,\Lambda_{ij}] = -(1+\delta_{kj})\Lambda_{ij} \]
\[[S_k,L_{ij}] = (-1)^{i+1}(1-3\delta_{kj})L_{ij},\quad
[S_k,\Lambda_{ij}] = (-1)^{i}(1-3\delta_{kj})\Lambda_{ij}\]
\[[S_k,V_j] = 0, \qquad [S_k,A_j] = 0\]
\end{cor}

Guided by all the explicit computations of the action on the isotypical
component $V = V_{-2}^{\oplus 6}$ made up to this point, we now define in terms of the natural geometric operators a set of Serre generators for the algebra $\gg$.
\begin{dfn}
\label{def:cartangen}
\[\begin{array}{lcl}
\mathbf{e}_1 = \frac{1}{4}[L_{20},A_1] & &\mathbf{f}_1 =
\frac{1}{4}[V_1,\Lambda_{20}]\\
~& ~ & \\
\mathbf{e}_2 = \frac{1}{4}[L_{22},A_0] & &\mathbf{f}_2 =
\frac{1}{4}[V_0,\Lambda_{22}]\\
~& ~ & \\
\mathbf{e}_3 = V_0 & &\mathbf{f}_3 =
A_0\\
~& ~ & \\
\mathbf{e}_4 = \frac{1}{4}[L_{12},A_0] & &\mathbf{f}_4 =
\frac{1}{4}[V_0,\Lambda_{12}]\\
~& ~ & \\
\mathbf{e}_5 = \frac{1}{4}[L_{10},A_1] & &\mathbf{f}_5 =
\frac{1}{4}[V_1,\Lambda_{10}]\\
\end{array}\]
~ \\
 Moreover, for all $i\in\{1,..,5\}$ we define $\mathbf{h}_i =
[\mathbf{e}_i,\mathbf{f}_i]$.
\end{dfn}
As the $\mathbf{e}_i$ have by construction associated matrix $e_{i+1}^i$ once restricted to $V$ and the $\mathbf{f}_i$ are their respective adjoints, one gets:
\begin{pro}
\label{pro:serrerelations} The operators
$\mathbf{e}_i,\mathbf{f}_j,\mathbf{h}_k$ satisfy the Serre relations
for $\slsc$ and the $\mathbf{h}_i$ span the Cartan subalgebra $\h$:
\[\begin{array}{l}\mathbf{h}_1 =
\frac{1}{2}\left(H_1 - H_2 - S_1 - S_2\right)\\
\mathbf{h}_2 = \frac{1}{2}\left(H_0 - H_1 + S_2\right)\\
\mathbf{h}_3 = \frac{1}{2}\left(-H_0 + H_1 + H_2\right)\\
\mathbf{h}_4 = \frac{1}{2}\left(H_0 - H_1 - S_2\right)\\
\mathbf{h}_5 = \frac{1}{2}\left(H_1 - H_2 + S_1 +
S_2\right)\end{array}\]
\end{pro}

It would be interesting as a last remark to identify in the list of quadratic invariants the geometric operators $L_{ij},\Lambda_{ij},V_j,A_j$, the algebra $\h$ and the $\so2$ generator $J$. To do this one could of course use the explicit matrices for the quadratic invariants once restricted to $V$, which are not difficult to compute. 
One can however get very quickly a qualitative picture by using the notion of multidegree which we now introduce.\\
The decomposition $T^*X = W_0 \oplus W_1\oplus
W_2$ induces naturally a multi-degree on $\bigwedge^*T^*_\C X$ with
values in $\Z^3$, which we indicate with $mdeg$. This follows from
the equation
\[\bigwedge^n T_\C^*X \cong \bigoplus_{p+q+r = n}
\bigwedge^p\left(W_0\otimes\C\right)\oplus
\bigwedge^q\left(W_1\otimes\C\right)\oplus
\bigwedge^r\left(W_2\otimes\C\right)\] We notice furthermore that
the (complexified) decomposition above is preserved by the operator
$J$, and therefore $mdeg$ commutes with the action of $\so2$. 
\begin{pro} \label{mdeg-gg} The operators
$L_j,V_j,\Lambda_j,A_j,H_j,S_j$ are $mdeg$-homogeneous, with
multi-degrees:
\[\begin{array}{lll} mdeg(L_0) = (0,1,1) & mdeg(L_1) = (1,0,1) & mdeg(L_2) = (1,1,0)\\
\\
mdeg(\Lambda_0) = (0,-1,-1) & mdeg(\Lambda_1) = (-1,0,-1) &
mdeg(\Lambda_2) = (-1,-1,0)\\ \\
mdeg(V_0) = (2,0,0) & mdeg(V_1) = (0,2,0) & mdeg(V_2) = (0,0,2)\\ \\
mdeg(A_0) = (-2,0,0) & mdeg(A_1) = (0,-2,0) & mdeg(A_2) = (0,0,-2)\\
\\
mdeg(H_0) = (0,0,0) & mdeg(H_1) = (0,0,0) & mdeg(H_2) = (0,0,0)\\ \\
mdeg(S_0) = (0,0,0) & mdeg(S_1) = (0,0,0) & mdeg(S_2) = (0,0,0)
\end{array}\]
\end{pro}
\dimo The values for mdeg for the $L_j$ and the $V_j$ follow
immediately from mdeg of the corresponding forms and the dual
(contraction) operators have opposite value of mdeg. The remaing
values can be computed using the additivity of mdeg with respect to
the bracket. \qed 
\begin{pro}
Let $\{j,k,l\} = \{0,1,2\}$. Then\\
\[\begin{array}{lcl}Span\left(L_{1j},L_{2j}\right) &  =  & Span\left( [E_{w_k},E_{\overline{w}_l}], [E_{w_l},E_{\overline{w}_k}]\right)\\
~\\
Span\left(\Lambda_{1j},\Lambda_{2j}\right) &  = & Span\left( [I_{w_k},I_{\overline{w}_l}], [I_{w_l},I_{\overline{w}_k}]\right)\\
~\\
Span\left(V_j\right) & = & Span\left([E_{w_j},E_{\overline{w}_j}]\right)\\
~\\
Span\left(A_j\right) & = & Span\left([I_{w_j},I_{\overline{w}_j}]\right)\\
~\\
\h\oplus Span\left(J\right) & = & 
{\displaystyle \bigoplus_{m=0}^2}~ Span\left([E_{w_m},I_{w_m}],[E_{\overline{w}_m},I_{\overline{w}_m}]\right)

\end{array}\]
\end{pro}
\dimo The $mdeg$ of the $L_{ij}$ is the same of the corresponding $L_j$, and similarly for their adjoints. The mdegs of the quadratic monomials are immediately computed as they are the sum of those of their components. For example, $mdeg(E_{w_0}) = mdeg(E_{\overline{w}_0}) = (1,0,0)$ , $mdeg(E_{w_1}) = mdeg(E_{\overline{w}_1}) = (0,1,0)$ and therefore 
$mdeg([E_{w_0},E_{\overline{w}_1}] = (1,1,0)$, equal to that of $L_{12}$ and $L_{22}$. \qed

\end{document}